\documentclass[reqno,10pt]{amsart}
\usepackage{epsfig,psfrag,verbatim,amssymb,amsfonts,amsmath,pstricks}
\newtheorem{theorem}{Theorem}
\newtheorem{lemma}[theorem]{Lemma}

\newcommand{\field}[1]{\mathbb{#1}}

\newcommand{\R}{\field{R}}

\newcommand{\Z}{\field{Z}}
\newcommand{\N}{\field{N}}

\newcommand{\PP}{\field{P}}

\newcommand{\EE}{\field{E}}

\newcommand{\F}{\mathcal{F}}

\newcommand{\si}{\sigma}

\newcommand{\om}{\omega}
\newcommand{\Om}{\Omega}
\newcommand{\eps}{\varepsilon}
\newcommand{\won}{{\boldsymbol 1}}

\newcounter{constante}
\newcommand{\con}[1]{
\immediate\write 1{\noexpand\newlabel{#1}{{\theconstante}{\theconstante}}}
                    c_{\theconstante}
                    \stepcounter{constante}
                   }
\newcommand{\abel}[1]{}

\begin{document}

\setcounter{page}{1}

\title[0-1 law for planar RWRE]
{The zero-one law for planar random walks in i.i.d.\ random environments revisited}
\thanks{\textit{2000 Mathematics Subject Classification.} Primary 60K37, secondary 60F20.}
\thanks{\textit{Key words:}\quad Random Environment, Random Walk, Zero-One Law, Transience
 }

\maketitle

\vspace*{-2mm}
\begin{center}

{\sc By  Martin P.W.\ Zerner
}
\end{center}\vspace*{5mm}

{\footnotesize {\sc Abstract}. 
In this note we present a simplified proof of the zero-one law by Merkl and Zerner (2001) for
 directional transience of random walks in i.i.d.\ random environments (RWRE) on $\Z^2$.
Also, we indicate how to construct a two-dimensional counterexample in a non-uniformly elliptic and stationary environment which has better ergodic properties than the example given by Merkl and Zerner.
}\vspace*{5mm}

\section{Introduction}
Let us first recall the model of random walks in  random environments (RWRE), see also \cite{zei} for a survey.
For $d\geq 1$, we denote by $\mathcal P$ the set of $2d$-dimensional 
probability vectors, and
set $\Omega=\mathcal P^{\Z^d}$. Any $\om \in \Omega$, written as 
$\om=((\om(x,x+e))_{|e|=1})_{x\in\Z^d}$, will be called an {\it environment}.
It is called
{\it elliptic} if $\om(x,x+e)>0$ 
for all $x,e\in \Z^d$ with $ |e|=1$
and
{\it uniformly elliptic} 
if there exists a so-called
{\it ellipticity constant}
 $\kappa>0$, such that $\om(x,x+e)>\kappa$ for 
all $x,e\in \Z^d$ with $ |e|=1$.
Endowing $\Om$ with the canonical product $\si$-algebra and a probability measure $\PP$ turns $\om$ into a collection of random $2d$-vectors, i.e.\ a {\it random environment}. The expectation corresponding to $\PP$ is denoted by $\EE$.

Given an environment $\om\in\Om$,
the values $\om(x,x+e)$ serve as transition probabilities for the $\Z^d$-valued Markov chain
$(X_n)_{n\geq 0}$, called {\it random walk in random environment} (RWRE). This 
process can be defined as the sequence of canonical projections on the sample space $(\Z^d)^\N$ endowed with the so-called {\it quenched measure} $P_{z,\om}$, which is defined for any starting point $z\in\Z^d$ and any environment $\om\in\Omega$ and characterized by
\begin{eqnarray*}
P_{z,\om}[X_0=z]&=&1\quad\mbox{and}\\
P_{z,\om}[X_{n+1}=X_n+e\ |\ X_0, X_1,\ldots,X_n]&=&\om(X_n,X_n+e)\quad P_{z,\om}-a.s.
\end{eqnarray*}
for all $e\in\Z^d$ with $|e|=1$ and all $n\geq 0$.
The so-called {\it annealed measures} $P_z$,\ $z\in\Z^d,$ are then defined as the semi-direct
products $P_z:=\PP\times P_{z,\om}$
on $\Omega\times(\Z^d)^\N$ by $P_z[\cdot]:=\EE[P_{z,\om}[\cdot]]$. The expectations corresponding to 
$P_{z,\om}$ and $P_z$ are denoted by $E_{z,\om}$ and $E_z$, respectively.

One of the major open questions in the study of  RWRE
concerns the so-called {\it 0-1 law}, which we shall describe now. For $\ell\in \R^d$, $\ell\neq 0$,
 define the event
\[A_{\ell}:=\left\{\lim_{n\to\infty}X_n\cdot \ell=\infty\right\}
\]
that the walker tends in a rough sense into direction $\ell$, which we call \textit{to the right}.
It has been 
known since the work of Kalikow \cite{Ka81},
that if the random vectors 
$\om(x,\cdot)$, $x\in\Z^d$, are
i.i.d.\ under $\PP$ and $\PP$-a.s.\ 
uniformly elliptic, then
\begin{equation}\label{kali}
 P_0\left[A_{\ell}\cup A_{-\ell}\right]\in \{0,1\}\,.\end{equation}
This was extended in \cite[Proposition 3]{ZeM01} to the 
elliptic i.i.d.\  case. We shall call (\ref{kali}) Kalikow's zero-one law.

The zero-one law for directional transience is the stronger statement that even  $P_0[A_{\ell}]\in\{0,1\}$.
Except for $d=1$, see e.g.\ \cite[Theorem 2.1.2]{zei},
it is only partially known under which conditions this statement holds.  For $d=2$ and $\om(x,\cdot)$, $x\in\Z^d$, being
i.i.d.\ under $\PP$, Kalikow \cite{Ka81}
presented it as a question;  that 
case was settled 
in the affirmative in \cite{ZeM01}, while the case $d\ge 3$ is still wide open.
\begin{theorem}\label{01}{\rm (}see \cite[Theorem 1]{ZeM01}{\rm )}
Let $d=2$, $\ell\in\R^2\backslash\{0\}$ and let $(\om(x,\cdot)),\ x\in\Z^2,$ be  i.i.d.\ and elliptic under $\PP$. Then
$P_0[A_{\ell}]\in\{0,1\}.$
\end{theorem} \abel{01}
 Let us briefly sketch the proof of 
 Theorem \ref{01} for $\ell=e_1$ as given in  \cite{ZeM01}, where $e_1$ is the first coordinate direction.
Assuming that the zero-one law does not hold, one considers two independent random walks in the same environment.
The first one starts at the origin, the second one at a point $(L,z_L)$ for some  $L>0$ large. The slab $\{(x_1,x_2)\in\Z^2\mid 0\le x_1\le L\}$ is then subdivided into three slabs of equal size. By adjusting $z_L$ and using $d=2$ one can then force the paths of the two walkers to intersect at some point $x$ in the middle slab with a positive probability, which is bounded away from 0 uniformly in $L$. In this step some technical result \cite[Lemma 7]{ZeM01} about sums of four independent random variables is used. Now consider a third random walker starting at $x$.
By Kalikow's zero-one law (\ref{kali}) it eventually needs to go either to the left or to the right. Since $x$ has been visited by the first walker which has already traveled a long distance $>L$ to the right and thus most probably will continue to go to the right, the third walker is also likely to go right. However, by the same argument the third walker should also follow the second walker to the left. This leads to the desired contradiction.

The main goal of the present paper is to give a simplified proof of Theorem \ref{01}. In Section 2 we are going to present a proof in which the slab between 0 and $(L,z_L)$ is divided into two slab only. This way the technical lemma \cite[Lemma 7]{ZeM01} is not needed anymore and general directions $\ell\notin\{e_1,e_2\}$ can be handled more easily.

The
same paper \cite{ZeM01} also provides an example
of a stationary
elliptic, ergodic  environment
with $P_0[A_{\ell}]\notin\{0,1\}.$ 
However, the environment in this example  has bad mixing properties. In fact, 
it is not even totally ergodic since it is not ergodic with respect to the subgroup of shifts in $2\Z^2$.
In Section 3 we shall sketch an alternative construction of  a counterexample which has better mixing properties than the one given in \cite{ZeM01}.

\section{A shorter proof of Theorem \ref{01}}
Without loss of generality we assume $\|\ell\|_2=1$.
For $u\in\R$ and $\diamond \in\{<,\le, >,\ge\}$ we consider the stopping times
\[T_{\diamond u}:=\inf\{n\ge 0\mid (X_n\cdot\ell) \diamond u\}.\]
By Kalikow's 0-1 law (\ref{kali}),
$P_0\left[A_{\ell}\cup A_{-\ell} \right]\in\{0,1\}$.
The case $P_0[A_{\ell}\cup A_{-\ell}]=0$ is trivial. So assume 
\begin{equation}\label{star}
P_0[A_{\ell}\cup A_{-\ell}]=1.
\end{equation}\abel{star}
For the proof of the theorem it suffices to show that\footnote{Here and in the following the sole purpose of the figures is to illustrate the term immediately preceeding the figure.
The proof is complete without the figures.}
\psset{unit=3mm}
\begin{equation}\label{kal}
0
=P_0[T_{<0}=\infty]\ P_0[T_{>0}=\infty] =
P\Bigg[
\begin{pspicture}(0.4,.7)(5,2)
\psline[linestyle=dotted](1.5,-.5)(1.5,2.5)
\psbezier[linewidth=0.2mm]{*->}%
(1.5,1)(2.5,4)(4,-2)(5,1)
\rput(.8,1){0}
\end{pspicture}
\Bigg]\ 
P\Bigg[
\begin{pspicture}(0.4,.7)(5,2)
\psline[linestyle=dotted](4,-.5)(4,2.5)
\psbezier[linewidth=0.2mm]{<-*}%
(.5,1)(1.5,4)(3,-2)(4,1)
\rput(4.7,1){0}
\end{pspicture}
\Bigg].
\end{equation}\abel{kal}
Indeed, that (\ref{kal}) is sufficient follows from \cite[Proposition 1.2 (1.16)]{SzZe99}. For completeness, we repeat the argument here.
If (\ref{kal}) holds then 
 either  $P_0$-a.s.\ $T_{<0}<\infty$ or  $P_0$-a.s.\ $T_{>0}<\infty$.
In the first case, due to translation invariance, $T_{<x}<\infty$ would  hold $P_x$-a.s.\ for all $x\in\Z^2$.
Hence, $\PP$-a.s.\ $P_{x,\om}[T_{<x}<\infty]=1$  for all $x\in\Z^2$. By the strong Markov property this implies
$P_0[A_{\ell}]=0$. Similarly, we get in the second case that $P_0[A_{-\ell}]=0$, which yields due to (\ref{star})
$P_0[A_{\ell}]=1-P_0[A_{-\ell}]=1$.

For the proof of (\ref{kal}) observe that 
\begin{equation}\label{near}
T_{\ge L},T_{\le -L}\ge L\quad\mbox{$P_0$-a.s.\ for all $L\ge 0$},
\end{equation}\abel{near} since the walk is moving between nearest neighbors.
Therefore, the right hand side of (\ref{kal}) is equal to
\begin{equation}\label{LL}
\lim_{L\to\infty} P_0[T_{\ge L}<T_{<0}]\ P_0[T_{\le -L}<T_{>0}]
 =\lim_{L\to\infty} 
P\Bigg[
\begin{pspicture}(0.2,.7)(5.8,2)
\psline[linestyle=dotted](1.5,-.5)(1.5,2.5)
\psline[linestyle=dotted](5,-.5)(5,2.5)
\psline[linewidth=0.1mm]{<-}(1.6,0)(3,0)
\rput(3.25,0){$L$}
\psline[linewidth=0.1mm]{->}(3.5,0)(4.9,0)
\psbezier[linewidth=0.2mm]{*->}%
(1.5,1)(2.5,1.5)(4,.5)(5,1)
\rput(.8,1){0}
\end{pspicture}
\Bigg]
P\Bigg[
\begin{pspicture}(-.3,.7)(5.3,2)
\psline[linestyle=dotted](4,-.5)(4,2.5)
\psline[linestyle=dotted](.5,-.5)(.5,2.5)
\psline[linewidth=0.1mm]{<-}(.6,0)(2,0)
\rput(2.25,0){$L$}
\psline[linewidth=0.1mm]{->}(2.5,0)(3.9,0)
\psbezier[linewidth=0.2mm]{<-*}%
(.5,1)(1.5,1.5)(3,.5)(4,1)
\rput(4.7,1){0}
\end{pspicture}
\Bigg].
\end{equation}\abel{LL}
Now fix  a unit vector $\ell^\perp\in\R^2$ which is perpendicular to $\ell$. Then for $L\in\N$ we choose 
$z_L\in\Z^2$ 
such that 
\begin{eqnarray}
&&x_L:=z_L\cdot \ell\ge 2L,\nonumber\\
&&\mbox{$z_L$ has a nearest neighbor $w_L\in\Z^2$ with $w_L\cdot\ell<2L$},\nonumber\\
\label{less}&&\mbox{$y_L:=z_L\cdot \ell^\perp$ is
a median of the distribution of $X_{T_{\ge 2L}}\cdot \ell^\perp$
under }\\ 
\nonumber
&&\mbox{\qquad $P_0[\ \cdot\mid T_{\ge 2L}<T_{<0}]$, i.e.\ }P_0[X_{T_{\ge 2L}}\cdot \ell^\perp\diamond y_L\mid T_{\ge 2L}<T_{<0}]\le 1/2\quad\mbox{for}\
\diamond\in\{<,>\}.
\end{eqnarray}
(If $\ell=e_1$ then we can choose $\ell^\perp=e_2$, $x_L=2L$ and $z_L=(2L,y_L)$.) 
In order to make the events in the two probabilities in (\ref{LL}) depend on disjoint and therefore independent parts of the environment we shift the starting point in the second factor in (\ref{LL}) from 0 to $z_L$. Thus, by translation invariance, we can rewrite (\ref{LL}) as 
\begin{equation}\label{LLL}
\lim_{L\to\infty} P_0[T_{\ge L}<T_{<0}]\ P_{z_L}[T_{\le L}<T_{>x_L}]
 =\lim_{L\to\infty} 
P\Bigg[
\begin{pspicture}(0.2,.7)(5.8,2)
\psline[linestyle=dotted](1.5,-.5)(1.5,2.5)
\psline[linestyle=dotted](5,-.5)(5,2.5)
\psline[linewidth=0.1mm]{<-}(1.6,0)(3,0)
\rput(3.25,0){$L$}
\psline[linewidth=0.1mm]{->}(3.5,0)(4.9,0)
\psbezier[linewidth=0.2mm]{*->}%
(1.5,1)(2.5,1.5)(4,.5)(5,1)
\rput(.8,1){0}
\end{pspicture}
\Bigg]
P\Bigg[
\begin{pspicture}(-.3,.7)(5.6,2)
\psline[linestyle=dotted](4,-.5)(4,2.5)
\psline[linestyle=dotted](.5,-.5)(.5,2.5)
\psline[linewidth=0.1mm]{<-}(.6,0)(2,0)
\rput(2.25,0){$L$}
\psline[linewidth=0.1mm]{->}(2.5,0)(3.9,0)
\psbezier[linewidth=0.2mm]{<-*}%
(.5,1)(1.5,1.5)(3,.5)(4,1)
\rput(5,1){$z_L$}
\end{pspicture}
\Bigg].
\end{equation}\abel{LLL}
To write the product of probabilities in (\ref{LLL}) as a single probability, we  introduce two independent random walks moving in the same environment, one starting at 0, the other starting at $z_L$. So for any $\om\in\Omega$ and $L\in\N$ let $P_{0,z_L,\om}$ be a probability measure on $(\Z^2)^\N\times
 (\Z^2)^\N$ such that the two canonical processes of projections $(X_n^1)_n$ and $(X_n^2)_n$ on this space are independent of each other and have distributions $P_{0,\om}$ and $P_{z_L,\om}$, respectively, and denote by 
$P_{0,z_L}$ the corresponding annealed measure. 
Stopping times referring to the walk $(X_n^i)_n$ will be marked with an upper index $i\ (i=1,2)$.
Then
(\ref{LLL}) is equal to
\begin{equation}\label{zusammen}
\lim_{L\to\infty}P_{0,z_L}[T_{\ge L}^1<T_{<0}^1, T_{\le L}^2<T_{>x_L}^2]
=\lim_{L\to\infty}
P\Bigg[
\begin{pspicture}(0.2,.7)(10,2)
\psline[linestyle=dotted](1.5,-.5)(1.5,2.5)
\psline[linestyle=dotted](5,.5)(5,2.5)
\psline[linestyle=dotted](8.5,-.5)(8.5,2.5)
\psline[linewidth=0.1mm]{<-}(1.6,0)(4.5,0)
\rput(5,0){$2L$}
\psline[linewidth=0.1mm]{->}(5.5,0)(8.4,0)
\psbezier[linewidth=0.2mm]{*->}%
(1.5,1)(2.5,1.5)(4,.5)(5,1)
\rput(.8,1){0}
\psbezier[linewidth=0.2mm]{<-*}%
(5,1.5)(6,3)(7.5,-.5)(8.5,1.5)
\rput(9.5,1.5){$z_L$}
\end{pspicture}
\Bigg].
\end{equation}\abel{zusammen}
After crossing the line $\{x\mid x\cdot\ell=L\}$ any walk must due to (\ref{star}) a.s.\ cross  the line $\{x\mid x\cdot\ell =0\}$ or the line $\{x\mid x\cdot\ell =2L\}$.
Consequently, (\ref{zusammen}) is less than or equal to
\begin{eqnarray}\label{affe}
\lefteqn{\liminf_{L\to\infty}P_0[T_{\ge L}<T_{<0}<\infty]+P_{z_L}[T_{\le L}<T_{>x_L}<\infty]+
P_{0,z_L}[T_{\ge 2L}^1<T_{<0}^1, T_{\le 0}^2<T_{>x_L}^2]}\\
&=& \liminf_{L\to\infty}P\Bigg[
\begin{pspicture}(0.2,.7)(5.8,2)
\psline[linestyle=dotted](1.5,-.5)(1.5,2.5)
\psline[linestyle=dotted](5,-.5)(5,2.5)
\psbezier[linewidth=0.2mm]{*->}%
(1.5,1)(7,-1)(7,3)(1,1.5)
\rput(.8,1){0}
\end{pspicture}
\Bigg] +
P\Bigg[
\begin{pspicture}(-.3,.7)(5.4,2)
\psline[linestyle=dotted](4,-.5)(4,2.5)
\psline[linestyle=dotted](.5,-.5)(.5,2.5)
\psbezier[linewidth=0.2mm]{<-*}%
(4.5,1.5)(-1.5,3)(-1.5,-1)(4,1)
\rput(5,1){$z_L$}
\end{pspicture}
\Bigg]+
P\Bigg[
\begin{pspicture}(0.2,.7)(9.5,2)
\psline[linestyle=dotted](1.5,-.5)(1.5,2.5)
\psline[linestyle=dotted](5,-.5)(5,2.5)
\psline[linestyle=dotted](8.5,-.5)(8.5,2.5)
\psbezier[linewidth=0.2mm]{*->}%
(1.5,1)(3.5,2.5)(6.5,-.5)(8.5,1)
\rput(.8,1){0}
\end{pspicture}
,
\begin{pspicture}(0.4,.7)(10,2)
\psline[linestyle=dotted](1.5,-.5)(1.5,2.5)
\psline[linestyle=dotted](5,-.5)(5,2.5)
\psline[linestyle=dotted](8.5,-.5)(8.5,2.5)
\psbezier[linewidth=0.2mm]{<-*}%
(1.5,1)(3.5,2.5)(6.5,-.5)(8.5,1)
\rput(9.5,1){$z_L$}
\end{pspicture}
\Bigg].
\nonumber
\end{eqnarray}
Due to (\ref{near})
the first term in (\ref{affe}) is $\le P_0[\exists n\ge L:\ |X_n\cdot\ell|\le 1]$.
The same holds for the second term, which is  
$\le P_{0}[T_{\le -L}<T_{>0}<\infty]$ due to translation invariance.
Therefore,
both terms vanish as $L\to\infty$ due to (\ref{star}).
Summarizing, we have shown
\begin{equation}\label{swr3}
P_0[T_{<0}=\infty]\ P_0[T_{>0}=\infty]\leq \liminf_{L\to\infty}P_{0,z_L}[T_{\ge 2L}^1<T_{<0}^1, T_{\le 0}^2<T_{>x_L}^2].
\end{equation}\abel{swr3}
Now consider the event on the right-hand side of (\ref{swr3}). There are two possibilities:
Either the paths of the two random walks cross each other (with probability, say, $C_L$) before  $T_{\ge 2L}^1$ and $T_{\le 0}^2$, respectively, or they avoid each other 
(with probability, say, $N_L$).
Therefore, defining the hitting time of $x\in\Z^2$ by $H(x):=\inf\{n\ge 0\mid X_n=x\}$,
 we can rewrite the probability on the right-hand side of (\ref{swr3}) 
as $C_L+N_L$, 
where
\[
C_L:= P_{0,z_L}[\exists x:\ H^1(x)\le T_{\ge 2L}^1<T_{<0}^1,\ H^2(x)\le T_{\le 0}^2<T_{>x_L}^2] 
=P\Bigg[
\begin{pspicture}(0.2,.7)(10,2)
\psline[linestyle=dotted](1.5,-.5)(1.5,2.5)
\psline[linestyle=dotted](5,-.5)(5,2.5)
\psline[linestyle=dotted](8.5,-.5)(8.5,2.5)
\psbezier[linewidth=0.2mm]{*->}%
(1.5,1)(3.2,-2.5)(4.9,4.5)(7.6,1)
\rput(.8,1){0}
\rput(8,1){\Bigg\{}
\rput(2,1){\Bigg\}}
\psbezier[linewidth=0.2mm]{<-*}%
(2.4,1)(3.5,0)(6.5,3)(8.5,1.5)
\rput(9.5,1.5){$z_L$}
\end{pspicture}
\Bigg]
\]
and
\[
N_L:= P_{0,z_L}[T_{\ge 2L}^1<T_{<0}^1,\ T_{\le 0}^2<T_{>x_L}^2, \{X_n^1\mid \ n\le T_{\ge 2L}^1\}\cap
\{X_n^2\mid \ n\le T_{\le 0}^2\}=\emptyset].
\]
Hence we get from (\ref{swr3})
\begin{equation}\label{frank}
P_0[T_{<0}=\infty]\ P_0[T_{>0}=\infty]\leq \limsup_{L\to\infty}C_L+\limsup_{L\to\infty}N_L.
\end{equation}\abel{frank}

To estimate $N_L$,  observe that  on the event in the definition of $N_L$ we have that 
$y_L-X_{T_{\ge 2L}^1}^1\cdot \ell^\perp$ and $X_{T_{\le0}^2}^2\cdot \ell^\perp$ are  either both strictly positive
or both strictly negative since the dimension is equal to two.
Indeed, otherwise the two paths $\{X_n^1\mid \ n\le T_{\ge 2L}^1\}$ and
$\{X_n^2\mid \ n\le T_{\le 0}^2\}$ would intersect each other, since the two diagonals of any  planar quadrangle intersect each other. 
 Therefore,
\begin{eqnarray*}
N_L&=& \sum_{s=\pm 1}
P_{0,z_L}\Big[T_{\ge 2L}^1<T_{<0}^1,\ T_{\le 0}^2<T_{>x_L}^2, \{X_n^1\mid  n\le T_{\ge 2L}^1\}\cap
\{X_n^2\mid  n\le T_{\le 0}^2\}=\emptyset,\\
&&\hspace*{19mm}s={\rm sign}\left(y_L-X_{T_{\ge 2L}^1}^1\cdot \ell^\perp\right)={\rm sign}\left(
X_{T_{\le0}^2}^2\cdot \ell^\perp\right)
\Big] \\
&=&
P\Bigg[
\begin{pspicture}(0.2,.7)(10,2)
\psline[linestyle=dotted](1.5,-.5)(1.5,2.5)
\psline[linestyle=dotted](5,-.5)(5,2.5)
\psline[linestyle=dotted](8.5,-.5)(8.5,2.5)%
\psbezier[linewidth=0.2mm]{*->}%
(1.5,1)(3.6,1.8) (5.7,1.4) (7.8,2.2)
\rput(.8,1){0}
\rput(8.1,2.2){\{}
\rput(1.9,.2){\big\}}
\psbezier[linewidth=0.2mm]{<-*}%
(2.2,.2)(4.3,1) (6.4,.7) (8.5,1.5)
\rput(9.5,1.5){$z_L$}
\end{pspicture}
\Bigg]
+
P\Bigg[
\begin{pspicture}(0.2,.7)(10,2)
\psline[linestyle=dotted](1.5,-.5)(1.5,2.5)
\psline[linestyle=dotted](5,-.5)(5,2.5)
\psline[linestyle=dotted](8.5,-.5)(8.5,2.5)
\psbezier[linewidth=0.2mm]{*->}%
(1.5,1)(3.6,1.2) (5.7,.2) (7.8,.4)
\rput(.8,1){0}
\rput(8.1,.4){\Big\{}
\rput(1.9,1.8){\big\}}
\psbezier[linewidth=0.2mm]{<-*}%
(2.2,1.8)(4.3,2.1) (6.4,1.2) (8.5,1.5)
\rput(9.5,1.5){$z_L$}
\end{pspicture}
\Bigg]
\end{eqnarray*}
Denoting by $\Pi_{L,s}$ $(L\in \N, s\in\{-1,0,+1\})$ the set of all finite nearest-neighbor paths that start at $z_L$ and leave the strip
$\{x\mid 0\le x\cdot\ell \le x_L\}$ on the opposite side through a vertex $x$ with ${\rm sign}\,(x\cdot \ell^\perp)=s$,  we rewrite $N_L$ as
\begin{eqnarray*}
N_L&=& \sum_{s=\pm 1}\sum_{\pi\in \Pi_{L,s}}P_{0,z_L}\Big[
T_{\ge 2L}^1<T_{<0}^1,\{X_n^1\mid  n\le T_{\ge 2L}^1\}\cap \pi=\emptyset,\ \mbox{$(X^2_n)_n$ follows $\pi$},\\
&& \hspace*{29mm} s={\rm sign}\left(y_L-X_{T_{\ge 2L}^1}^1\cdot \ell^\perp\right)\Big].
\end{eqnarray*}
Using the disjointness of the paths we get by independence in the environment,
\begin{eqnarray*}
N_L&=&\sum_{s=\pm 1}\sum_{\pi\in \Pi_{L,s}}
P_0\left[T_{\ge 2L}<T_{<0}, \{X_n^1\mid  n\le T_{\ge 2L}^1\}\cap \pi=\emptyset,
s={\rm sign}\left(y_L-X_{T_{\ge 2L}}\cdot \ell^\perp\right)\right]\\
&& \hspace*{18mm}P_{z_L}\left[\mbox{$(X_n)_n$ follows $\pi$}\right]\\
&\le &\sum_{s=\pm 1}\sum_{\pi\in \Pi_{L,s}} 
P_0\left[T_{\ge 2L}<T_{<0},  s={\rm sign}\left(y_L-X_{T_{\ge 2L}}\cdot \ell^\perp\right)\right]\
 P_{z_L}\left[\mbox{$(X_n)_n$ follows $\pi$}\right]\\
&=&\sum_{s=\pm 1}P_0\left[T_{\ge 2L}<T_{<0}, s={\rm sign}\left(y_L-X_{T_{\ge 2L}}\cdot \ell^\perp\right)\right]\\
&& \hspace*{8mm}P_{z_L}\left[T_{\le 0}<T_{>x_L},s={\rm sign}\left(X_{T_{\le 0}}\cdot \ell^\perp\right)\right]
\end{eqnarray*}
\begin{eqnarray*}
&=&P\Bigg[
\begin{pspicture}(0.2,.7)(10,2)
\psline[linestyle=dotted](1.5,-.5)(1.5,2.5)
\psline[linestyle=dotted](5,-.5)(5,2.5)
\psline[linestyle=dotted](8.5,-.5)(8.5,2.5)
\psbezier[linewidth=0.2mm]{*->}%
(1.5,1)(3.6,3.4) (5.8,-.2) (7.9,2.2)
\rput(.8,1){0}
\rput(8.1,2.2){\{}
\psline[linewidth=0.2mm]{-*}%
(8.5,1.5)(8.5,1.5)
\rput(9.5,1.5){$z_L$}
\end{pspicture}
\Bigg]
P\Bigg[
\begin{pspicture}(0.2,.7)(10,2)
\psline[linestyle=dotted](1.5,-.5)(1.5,2.5)
\psline[linestyle=dotted](5,-.5)(5,2.5)
\psline[linestyle=dotted](8.5,-.5)(8.5,2.5)
\psline[linewidth=0.2mm]{*-}%
(1.5,1)(1.5,1)
\rput(.8,1){0}
\rput(1.9,.2){\big\}}
\psbezier[linewidth=0.2mm]{<-*}%
(2.2,.2)(4.3,-1.4) (6.4,3.1) (8.5,1.5)
\rput(9.5,1.5){$z_L$}
\end{pspicture}
\Bigg]\\
&&
+\ 
P\Bigg[
\begin{pspicture}(0.2,.7)(10,2)
\psline[linestyle=dotted](1.5,-.5)(1.5,2.5)
\psline[linestyle=dotted](5,-.5)(5,2.5)
\psline[linestyle=dotted](8.5,-.5)(8.5,2.5)
\psbezier[linewidth=0.2mm]{*->}%
(1.5,1)(3.6,2.8)(5.6,-1.4) (7.7,.4)
\rput(.8,1){0}
\rput(8.1,.4){\Big\{}
\psline[linewidth=0.2mm]{-*}%
(8.5,1.5)(8.5,1.5)
\rput(9.5,1.5){$z_L$}
\end{pspicture}
\Bigg]
P\Bigg[
\begin{pspicture}(0.2,.7)(10,2)
\psline[linestyle=dotted](1.5,-.5)(1.5,2.5)
\psline[linestyle=dotted](5,-.5)(5,2.5)
\psline[linestyle=dotted](8.5,-.5)(8.5,2.5)
\psline[linewidth=0.2mm]{*-}%
(1.5,1)(1.5,1)
\rput(.8,1){0}
\rput(1.9,1.8){\big\}}
\psbezier[linewidth=0.2mm]{<-*}%
(2.2,1.8)(4.2,-.3)(6.4,3.6 )(8.5,1.5)
\rput(9.5,1.5){$z_L$}
\end{pspicture}
\Bigg]\\
&\stackrel{(\ref{less})}{\le}& \frac{1}{2} P_0\left[T_{\ge 2L}<T_{<0}\right]\sum_{s=\pm 1} 
P_{z_L}\left[T_{\le 0}<T_{>x_L},s={\rm sign}\left(X_{T_{\le 0}}\cdot \ell^\perp\right)\right]\\
&=& \frac{1}{2} 
P\Bigg[
\begin{pspicture}(0.2,.7)(9.5,2)
\psline[linestyle=dotted](1.5,-.5)(1.5,2.5)
\psline[linestyle=dotted](5,-.5)(5,2.5)
\psline[linestyle=dotted](8.5,-.5)(8.5,2.5)
\psbezier[linewidth=0.2mm]{*->}%
(1.5,1)(3.5,2.5)(6.5,-.5)(8.5,1)
\rput(.8,1){0}
\end{pspicture}
\Bigg]\ 
\left(
P\Bigg[
\begin{pspicture}(0.2,.7)(10,2)
\psline[linestyle=dotted](1.5,-.5)(1.5,2.5)
\psline[linestyle=dotted](5,-.5)(5,2.5)
\psline[linestyle=dotted](8.5,-.5)(8.5,2.5)
\psline[linewidth=0.2mm]{*-}%
(1.5,1)(1.5,1)
\rput(.8,1){0}
\rput(1.9,.2){\big\}}
\psbezier[linewidth=0.2mm]{<-*}%
(2.2,.2)(4.3,-1.4) (6.4,3.1) (8.5,1.5)
\rput(9.5,1.5){$z_L$}
\end{pspicture}
\Bigg]
+ 
P\Bigg[
\begin{pspicture}(0.2,.7)(10,2)
\psline[linestyle=dotted](1.5,-.5)(1.5,2.5)
\psline[linestyle=dotted](5,-.5)(5,2.5)
\psline[linestyle=dotted](8.5,-.5)(8.5,2.5)
\psline[linewidth=0.2mm]{*-}%
(1.5,1)(1.5,1)
\rput(.8,1){0}
\rput(1.9,1.8){\big\}}
\psbezier[linewidth=0.2mm]{<-*}%
(2.2,1.8)(4.2,-.3)(6.4,3.6 )(8.5,1.5)
\rput(9.5,1.5){$z_L$}
\end{pspicture}
\Bigg]
\right)\\
&\le& \frac{1}{2}P_0\left[T_{\ge 2L}<T_{<0}\right]\ P_{z_L}\left[T_{\le 0}<T_{>x_L}\right]\\
&=&\frac{1}{2}
P\Bigg[
\begin{pspicture}(0.2,.7)(9.5,2)
\psline[linestyle=dotted](1.5,-.5)(1.5,2.5)
\psline[linestyle=dotted](5,-.5)(5,2.5)
\psline[linestyle=dotted](8.5,-.5)(8.5,2.5)
\psbezier[linewidth=0.2mm]{*->}%
(1.5,1)(3.5,2.5)(6.5,-.5)(8.5,1)
\rput(.8,1){0}
\end{pspicture}
\Bigg]\ 
P\Bigg[
\begin{pspicture}(0.4,.7)(10,2)
\psline[linestyle=dotted](1.5,-.5)(1.5,2.5)
\psline[linestyle=dotted](5,-.5)(5,2.5)
\psline[linestyle=dotted](8.5,-.5)(8.5,2.5)
\psbezier[linewidth=0.2mm]{<-*}%
(1.5,1)(3.5,2.5)(6.5,-.5)(8.5,1)
\rput(9.5,1){$z_L$}
\end{pspicture}
\Bigg]\\
&\le &\frac{1}{2}P_0\left[T_{\ge 2L}<T_{<0}\right]\ P_{0}\left[T_{\le -2L}<T_{>0}\right]\\
&=&\frac{1}{2}
P\Bigg[
\begin{pspicture}(0.2,.7)(9.5,2)
\psline[linestyle=dotted](1.5,-.5)(1.5,2.5)
\psline[linestyle=dotted](5,-.5)(5,2.5)
\psline[linestyle=dotted](8.5,-.5)(8.5,2.5)
\psbezier[linewidth=0.2mm]{*->}%
(1.5,1)(3.5,2.5)(6.5,-.5)(8.5,1)
\rput(.8,1){0}
\end{pspicture}
\Bigg]\ 
P\Bigg[
\begin{pspicture}(0.4,.7)(10,2)
\psline[linestyle=dotted](1.5,-.5)(1.5,2.5)
\psline[linestyle=dotted](5,-.5)(5,2.5)
\psline[linestyle=dotted](8.5,-.5)(8.5,2.5)
\psbezier[linewidth=0.2mm]{<-*}%
(1.5,1)(3.5,2.5)(6.5,-.5)(8.5,1)
\rput(9.5,1){$0$}
\end{pspicture}
\Bigg]\\
&\stackrel{(\ref{near})}{\le}& \frac{1}{2}P_0[2L<T_{<0}]\ P_0[2L<T_{>0}]\\
&\begin{array}[t]{c}\longrightarrow\vspace*{-2.5mm}\\
{\scriptscriptstyle L\to\infty}\end{array}
&\frac{1}{2}P_0[T_{<0}=\infty]\ P_0[T_{>0}=\infty]\ =\ 
\frac{1}{2}P\Bigg[
\begin{pspicture}(0.4,.7)(5,2)
\psline[linestyle=dotted](1.5,-.5)(1.5,2.5)
\psbezier[linewidth=0.2mm]{*->}%
(1.5,1)(2.5,4)(4,-2)(5,1)
\rput(.8,1){0}
\end{pspicture}
\Bigg]\ 
P\Bigg[
\begin{pspicture}(0.4,.7)(5,2)
\psline[linestyle=dotted](4,-.5)(4,2.5)
\psbezier[linewidth=0.2mm]{<-*}%
(.5,1)(1.5,4)(3,-2)(4,1)
\rput(4.7,1){0}
\end{pspicture}
\Bigg].
\end{eqnarray*}
Hence, due to (\ref{frank}),
\[\frac{1}{2}P_0[T_{<0}=\infty]\ P_0[T_{>0}=\infty]\leq \limsup_{L\to\infty}C_L.\]
For the proof of (\ref{kal}) it therefore suffices to show
\begin{equation}\label{doef}
\lim_{L\to\infty} C_L=0.
\end{equation}\abel{doef}
By considering the possible locations of the intersections of the two paths,
we estimate $C_L$ by 
\begin{eqnarray*}
C_L&\le& C_0^L+C_L^{x_L}
\ =\ P\Bigg[
\begin{pspicture}(0.2,.7)(10,2)
\psline[linestyle=dotted](1.5,-.5)(1.5,2.5)
\psline[linestyle=dotted](5,-.5)(5,2.5)
\psline[linestyle=dotted](8.5,-.5)(8.5,2.5)
\psbezier[linewidth=0.2mm]{*->}%
(1.5,1)(3.5,2.5)(6.5,-1)(7.6,1)
\psline[linewidth=0.2mm]{*-}(4,1.15)(4,1.15)
\rput(4,1.8){$x$}
\rput(.8,1){0}
\rput(8,1){\Bigg\{}
\rput(2,1){\Bigg\}}
\psbezier[linewidth=0.2mm]{<-*}%
(2.4,1)(3.5,0)(6.5,4)(8.5,1.5)
\rput(9.5,1.5){$z_L$}
\end{pspicture}
\Bigg] + P\Bigg[
\begin{pspicture}(0.2,.7)(10,2)
\psline[linestyle=dotted](1.5,-.5)(1.5,2.5)
\psline[linestyle=dotted](5,-.5)(5,2.5)
\psline[linestyle=dotted](8.5,-.5)(8.5,2.5)
\psbezier[linewidth=0.2mm]{*->}%
(1.5,1)(3.2,2)(4.9,3)(7.6,1)
\psline[linewidth=0.2mm]{*-}(6.25,1.83)(6.25,1.83)
\rput(6.25,1.2){$x$}
\rput(.8,1){0}
\rput(8,1){\Bigg\{}
\rput(2,1){\Bigg\}}
\psbezier[linewidth=0.2mm]{<-*}%
(2.4,1)(3.5,0)(6.5,3)(8.5,1.5)
\rput(9.5,1.5){$z_L$}
\end{pspicture}
\Bigg],
\end{eqnarray*} 
where
\[C_a^b:= P_{0,z_L}[\exists x:\ a\le x\cdot \ell\le b,\ H^1(x)\le T_{\ge 2L}^1<T_{<0}^1,\ H^2(x)\le T_{\le 0}^2<T_{>x_L}^2].\] 
Due to symmetry and translation invariance it suffices to show for the proof of (\ref{doef}) that $C_0^L\to\infty$.
To this end let $\eps>0$ and set $r(x,\om):=P_{x,\om}[A_{\ell}]$. Then 
\begin{eqnarray}
C_0^L&\le& C_{0,1}^L+ C_{0,2}^L,\quad\mbox{ where}\label{ell}\\ \nonumber
C_{0,1}^L&:=& P_{0}[\exists x:\ r(x,\om)\le \eps,\ H(x)\le T_{\ge L}<\infty]=
P\Bigg[
\begin{pspicture}(0.2,.7)(5.8,2)
\psline[linestyle=dotted](1.5,-.5)(1.5,2.5)
\psline[linestyle=dotted](5,-.5)(5,2.5)
\psbezier[linewidth=0.2mm]{*->}%
(1.5,1)(2.5,2.5)(4,-.5)(5,1)
\rput(.8,1){0}
\psline[linewidth=0.2mm]{*-}(2.5,1.45)(2.5,1.45)
\rput(2.5,.8){$x$}
\end{pspicture},
r(x,\om)\le \eps
\Bigg]
\quad\mbox{and}\\
C_{0,2}^L&:=& P_{z_L}[\exists x:\ x\cdot\ell\le L,\ r(x,\om)\ge \eps,\ H(x)<\infty]\ =\
P\Bigg[
\begin{pspicture}(1.7,.7)(10,2)
\psline[linestyle=dotted](5,-.5)(5,2.5)
\psline[linestyle=dotted](8.5,-.5)(8.5,2.5)
\psline[linewidth=0.2mm]{*-}(4,1.15)(4,1.15)
\rput(4,1.8){$x$}
\psbezier[linewidth=0.2mm]{<-*}%
(2.5,1)(3.5,0)(6.5,4)(8.5,1.5)
\rput(9.5,1.5){$z_L$}
\end{pspicture}, r(x,\om)\ge \eps
\Bigg].\nonumber
\end{eqnarray}
In order to bound $C_{0,1}^L$, consider
$\si:=\inf\{n\ge 0\mid r(X_n,\om)\le \eps\}.$ Note that $\si$ is a stopping time w.r.t.\ the filtration $(\F_n)_{n\ge 0}$, where $\F_n$ is the $\si$-field generated by $X_0,\ldots,X_n$ and $\om$.
 Therefore, by the strong Markov property,
\[
 C_{0,1}^L=P_0[\si\le T_{\ge L}<\infty]=E_0\left[P_{X_{\si},\om}[T_{\ge L}<\infty], 
\si\le T_{\ge L}, \si<\infty\right].
\]
 Since for all $x\in\Z^2$ and almost all $\om$, $P_{x,\om}$-a.s.\ $\{T_{\ge L}<\infty\}\searrow A_{\ell}$ as 
$L\to\infty$ due to (\ref{star}), we get 
by dominated convergence
\begin{equation}
\label{melli}
\lim_{L\to\infty} C_{0,1}^L=E_0\left[P_{X_{\si,\om}}[A_{\ell}], \si<\infty
\right]
 \ =\ 
E_0\left[r\left(X_{\si},\om\right),\ \si<\infty
\right]\ \le\ \eps
\end{equation}\abel{melli}
by definition of $\si$.
Now consider $C_{0,2}^L$. By translation invariance
\begin{eqnarray}\nonumber
C_{0,2}^L&\le &P_{0}[\exists x:\ x\cdot\ell\le -L,\ r(x,\om)\ge \eps, H(x)<\infty]\\
&\stackrel{(\ref{star}),(\ref{near})}{\le}&
P_0[T_{\le -L}<\infty, A_{\ell}]+P_0[\exists n\ge L: r(X_n,\om)\ge \eps, A_{-\ell}].\label{term}
\end{eqnarray} 
Obviously, the first term in (\ref{term}) goes to 0 as $L\to\infty$. The same holds for the second term
since due to the martingale convergence theorem, $P_0$-a.s.\ $r(X_n,\om)=P_0[A_{\ell}\mid \F_n]\to\won_{A_{\ell}}$ as $n\to\infty$, cf.\
\cite[Lemma 5]{ZeM01}.
Together with (\ref{ell}) and (\ref{melli}) this yields $\limsup_{L}C_0^L\le \eps$.
Letting $\eps\searrow 0$ gives $\lim_L C_0^L=0$. This finishes the proof of (\ref{doef}).
\hfill$\Box$\vspace{2mm}

\section{A stationary and totally ergodic counterexample}
The stationary and ergodic environment constructed in \cite{ZeM01} is based on two disjoint trees which together span $\Z^2$. The branches of these trees are paths of coalescing random walks which for one tree go either up or right and for the other tree go either left or down. 
In order to allocate enough space for both trees some periodicity was introduced which destroyed total ergodicity of the environment. 

\begin{figure}[t]
\hspace*{-10mm}\epsfig{file=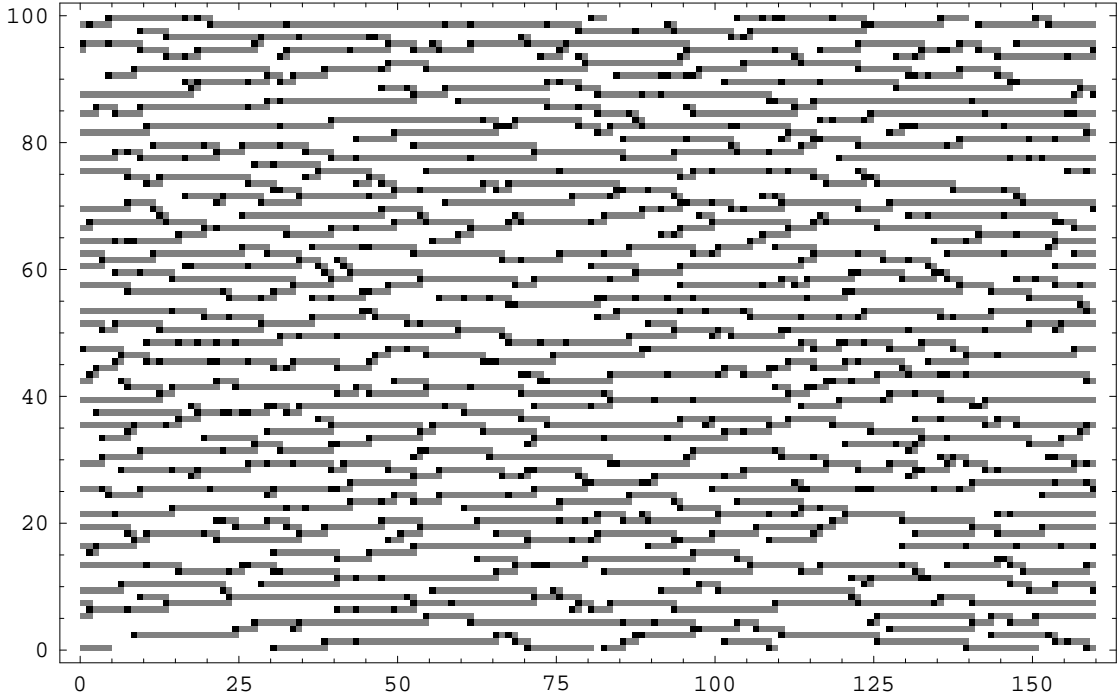,height=9.5cm, angle=0}
\caption{The black points constitute some in a certain way thinned discrete Poisson process on a torus. From each black point a gray line grew to the right until its tip became a $\ell^1$-neighbor of a black point.}
\label{paths}
\end{figure}

In this section we shall sketch an alternative construction which gives a totally ergodic environment. 
It has been inspired by the Poisson tree considered in \cite[Section 3]{FLT04}.
The main difference to the tree in \cite{FLT04} is that the nodes and leaves of our tree do not form a Poisson process but a point  process which has been obtained from a Poisson process by a local thinning procedure as follows:
We first color the vertices of $\Z^2$  independently black
with some fixed probability $0<p<1$ (in Figure \ref{paths}, $p=1/7$) and white otherwise. Then those black points $x\in\Z^2$ for which 
the set $x+\{\pm e_2, \pm 2e_2, -e_1\pm 2e_2\}$ contains another black point are removed
simultaneously, i.e.\ 
painted white again. The remaining set of black points is called  $B\subseteq\Z^2$. Obviously, the random variables $\mathbf 1_{x\in B},\ x\in\Z^2,$ are only finite range dependent.

Now each black point grows a gray line to the right until the line's tip reaches a neighbor of another black point,  see also Figure \ref{paths}.
The set of the gray points obtained this way is called $G\subseteq\Z^2$.
More formally, for $x\in\Z^2$ let 
\begin{eqnarray*}
g(x)&:=&\inf\left\{n\ge 0\mid x+ne_1\in (B+\{e_2,-e_2,-e_1\})\right\}
\quad\mbox{and set}\\
G&:=&\{x+ke_1\mid x\in B, 1\le k\le g(x)\}.
\end{eqnarray*}
Note that almost surely all $g(x),\ x\in\Z^2,$ are finite.
 Now consider the set $T:=B\cup G$. 
\abel{paths}
\begin{lemma}\label{T}
If $x\in T$ then $x+e_1\in T$ or $\{x+e_2,x+e_1+e_2\}\subseteq T$ or $\{x-e_2,x+e_1-e_2\}\subseteq T$.
Similarly, if $x\in T^c$ then $x-e_1\in T^c$ or  $\{x+e_2,x-e_1+e_2\}\subseteq T^c$ or $\{x-e_2,x-e_1-e_2\}\subseteq T^c$.
Moreover, $T^c\ne\emptyset$.
\end{lemma}\abel{T}
\begin{proof}
First note that
\begin{equation}\label{impl1}
x\in B\Rightarrow x+e_1\in T
\end{equation}\abel{impl1}
since either $x+e_1$ is black or it will be painted gray as the right neighbor of a black point.

For the first statement of the lemma let $x\in T$ and assume $x+e_1\in T^c$.  Then by (\ref{impl1}), $x\in T\backslash B=G$.
Hence it suffices to show that 
\begin{equation}\label{impl2}
x\in G, x+e_1\in T^c \Rightarrow ((x+e_2\in B\wedge x+e_1+e_2\in T)\vee (x-e_2\in B\wedge x+e_1-e_2\in T)),
\end{equation}\abel{impl2}
where $\wedge$  and $\vee$ denote ``and" and  ``or".
Since $x$ is gray but $x+e_1$ is white, the gray line to which $x$ belongs must have stopped growing in $x$.
This means that one of the neighbors $x+e_1, x+e_2$ or $x-e_2$ must be black. Hence, since $x+e_1$ is white 
$x+e_2$ or $x-e_2$ must be black.
By construction of $B$, only one of them can be black.  By symmetry we may assume
$x+e_2\in B$.  Then, due to (\ref{impl1}), $x+e_2+e_1\in T$.  Thus (\ref{impl2}) has been shown and the first statement follows.

For the second statement of the lemma let $x\in T^c$ and assume $x-e_1\in T$. Then applying (\ref{impl1}) to $x-e_1$ instead of $x$ yields $x-e_1\in T\backslash B=G$. 
Consequently, an application of (\ref{impl2}) to $x-e_1$ instead of $x$ yields without loss of generality (due to symmetry) that
$x-e_1+e_2\in B$ and $x+e_2\in T$. Now it suffices to show that $x-e_2, x-e_1-e_2\in T^c$. 
Neither of these points can be black by construction of $B$ since $x-e_1+e_2$ is already black. So it suffices to show that neither of them is gray. This is done by contradiction.
Assume that one of them is gray.
If $x-e_2$ were gray then $x-e_1-e_2$ would have to be black or gray. By construction of $B$, $x-e_1-e_2$  cannot be black since $x-e_1+e_2$ is already black. Hence  $x-e_1-e_2$  would have to be gray, too.
So we may assume that  $x-e_1-e_2$ is gray. By construction of the gray lines,  there is some $k\ge 2$ such that $x-ke_1-e_2$ is black and all the points $x-ie_1-e_2$\ $(1\le i<k)$  in between are gray.
 Now recall that $x-e_1$ is gray, too. Hence by the same argument,   there is some $m\ge 2$ such that $x-me_1$ is black and all the points $x-ie_1$\ $(1\le i<m)$  in between are gray. By construction of $B$,
$k$ and $m$ cannot be equal, since this would give two black points, $x-ke_1-e_2$ and $x-ke_1$, on top of each other.
 So assume $2\le k<m$. The case $2\le m<k$ is treated similarly. Then the gray line starting at the black point $x-me_2$ passes through the neighbor $x-ke_2$ of the black vertex  $x-ke_1-e_2$. Hence, it has to stop there, i.e.
the next point $x-(k-1)e_1-e_2$ cannot be gray, which it is. This gives the desired contradiction and proves the second statement.

For the last statement $T^c\ne\emptyset$ we show that $x\in B$ implies $x+e_1+e_2\in T^c$ or
 $x+e_1-e_2\in T^c$. Firstly, by construction of $B$,  not both  $x+e_1+e_2$ and $x+e_1-e_2$ can be black.
Secondly, none of them can be gray. Indeed, assume that for example $x+e_1+e_2$ were gray. As above, this would imply that its left neighbor $x+e_2$ would be black or gray, too.  However,  by construction of $B$, $x+e_2$  cannot be black since $x$ is already black. Hence  $x+e_2$ and $x+e_1+e_2$ would belong to the same gray line starting at some black point $x-ke_1+e_2$ with $k\ge 1$. However, this line would have to stop in  $x+e_2$ and not extend to $x+e_1+e_2$ since $x+e_2$ is a neighbor of the black point $x$, which would give a contradiction.
\end{proof}

Now we define for  any  vertex $x\in\Z^2$ its ancestor $a(x)$ as follows: For $x\in T$ we set
\begin{eqnarray*}
a(x):=\left\{\begin{array}{cl}
x+e_1&\mbox{if $x+e_1\in T$,}\\
x+e_2&\mbox{else if $x+e_2,x+e_1+e_2\in T$,}\\
x-e_2&\mbox{else if $x-e_2,x+e_1-e_2 \in T$}
\end{array}\right.
\end{eqnarray*}
and for $x\in T^c$ we define
\begin{eqnarray*}
a(x):=\left\{\begin{array}{cl}
x-e_1&\mbox{if $x-e_1\in T^c$,}\\
x+e_2&\mbox{else if $x+e_2,x-e_1+e_2\in T^c$,}\\
x-e_2&\mbox{else if $x-e_2,x-e_1-e_2 \in T^c$.}
\end{array}\right.
\end{eqnarray*}
Due to Lemma \ref{T} the function $a:\Z^2\to\Z^2$ is well defined and determines two disjoint infinite trees with sets of vertices $T$ and $T^c$, respectively. (The thinning of the Poisson process at the beginning was necessary to prevent the black and gray tree to disconnect the white complement into finite pieces, possibly leaving it without an infinite component.) Moreover, 
\begin{equation}\label{doppel}
\begin{array}{ll}
(a^n(x)\cdot e_1)_{n\ge 0}\quad\mbox{is monotone increasing for $x\in T$ and decreasing for $x\in T^c$ with}\\
a^2(x)\cdot e_1\ge x\cdot e_1+1\quad \mbox{for $x\in T$ and}\quad  a^2(x)\cdot e_1\le x\cdot e_1-1\quad \mbox{for $x\in T^c$.}
\end{array}
\end{equation}\abel{doppel}

It can be shown, cf.\ \cite[Theorem 3.1(d)]{FLT04}, that all the branches in the tree on $T$ are almost surely finite, i.e.\ 
 the length $h(x):=\sup\{n\ge 0\mid\exists y\ a^n(y)=x\}$ of the longest line of descendants of $x$ in $T$ is almost surely finite for all $x\in T$.
Moreover, it can be shown, cf.\ \cite[Theorem 3.1(b)]{FLT04}, that the tree on $T$ is connected. This implies that all the branches of the white tree on $T^c$ have finite height $h(x)$ as well.
The rest of the construction is the same as in \cite[pp.\ 1730, 1732]{ZeM01}: We define the environment in terms of the ancestor function $a$ for $x,y\in\Z^2$ with $|x-y|=1$ by
\[\om(x,y)=\left\{\begin{array}{ll}
 {\displaystyle 1-3/(h^2(x)+4)}&\mbox{if $y=a(x)$}\\                 
 {\displaystyle    1/(h^2(x)+4)} &\mbox{else.}           \end{array}\right.
\]
By Borel Cantelli there is a positive constant $c$ such that $P_{x,\om}[\forall n\ X_n=a^n(x)]>c$ for all $x\in\Z^2$ and almost all $\om$.
Consequently, due to (\ref{doppel}), and since neither tree is empty,
\[
P_0\left[\liminf_{n\to\infty}\frac{X_n\cdot e_1}{n}\ge \frac{1}{2}\right]>0\quad\mbox{and}\quad
P_0\left[\liminf_{n\to\infty}\frac{X_n\cdot(-e_1)}{n}\ge \frac{1}{2}\right]>0
\]
and in particular $P_0[A_{e_1}]\notin\{0,1\}.$
Since $(\om(x,\cdot))_{x\in\Z^2}$ 
has been obtained from $B$ by the application of a deterministic function which commutes with all spatial shifts in $\Z^2$
and since $B$ itself is stationary and totally ergodic with respect to all shifts, 
$(\om(x,\cdot))_{x\in\Z^2}$ is stationary and totally ergodic as well.
We refrained from 
investigating the mixing properties of this environment, which has been done for a similar counterexample for $d\ge 3$ in  \cite{BZZ}.\vspace*{2mm}

{\bf Open problems:} The gap between positive and negative results concerning the directional zero-one law in $d=2$ could be narrowed by answering one of the following questions: (1) Is there  
a  stationary and ergodic counterexample to the directional zero-one law for $d=2$, which is uniformly elliptic?
For $d\ge 3$  there is such a counterexample, even in a polynomially mixing environment, see \cite{BZZ}. (2)
Can  the directional zero-one law for $d=2$ be extended to stationary, ergodic and uniformly elliptic environments which have weaker independence properties than finite range dependence?
\bibliographystyle{amsalpha}
\vspace*{5mm}
{\sc\small 
Mathematisches Institut\\
Universit\"at T\"ubingen\\
Auf der Morgenstelle 10\\
72076 T\"ubingen, Germany\\
E-Mail: {\rm martin.zerner@uni-tuebingen.de} }
\end{document}